\documentclass[a4paper,12pt]{amsart}
\usepackage{amssymb}
\usepackage{enumerate}
\usepackage{mathrsfs}
\usepackage{graphicx}
\usepackage{booktabs}
\usepackage{stmaryrd}
\usepackage[dvips]{hyperref}
\usepackage{graphpap}
\usepackage{defs-eng}

\title[Combinatorial Line Arrangements]{Invariants of Combinatorial Line Arrangements 
and Rybnikov's Example}
\author[E. Artal]
{Enrique ARTAL BARTOLO}
\address{Departamento de Matem\'aticas\\
Campus Plaza de San Francisco s/n\\
E-50009 Zaragoza SPAIN}
\email{artal@unizar.es}

\author[J. Carmona]{Jorge CARMONA RUBER}
\address
{Departamento de Sistemas inform\'aticos y programaci\'on\\
Universidad Complu\-tense\\
Ciudad Universitaria s/n\\
E-28040 Madrid SPAIN}
\email{jcarmona@sip.ucm.es}

\author[J.I. Cogolludo]{Jos\'e Ignacio COGOLLUDO AGUST\'IN}
\address{Departamento de Matem\'aticas\\
Campus Plaza de San Francisco s/n\\
E-50009 Zaragoza SPAIN}
\email{jicogo@unizar.es}

\author[M.A. Marco]{Miguel \'Angel MARCO BUZUN\'ARIZ}
\address{Departamento de Matem\'aticas\\
Campus Plaza de San Francisco s/n\\
E-50009 Zaragoza SPAIN}
\email{mmarco@unizar.es}

\thanks{Second author is partially
supported by BFM2001-1488-C02-01; others authors are partially supported by
BFM2001-1488-C02-02.}

\begin{document}

\begin{abstract}
Following the general strategy proposed by G.Rybnikov, we present a
proof of his well-known result, that is, the existence of two arrangements of lines
having the same combinatorial type, but non-isomorphic fundamental groups. To do so,
the Alexander Invariant and certain invariants of combinatorial line arrangements
are presented and developed for combinatorics with only double and triple points.
This is part of a more general project to better understand the relationship 
between topology and combinatorics of line arrangements.
\end{abstract}

\maketitle

One of the main subjects in the theory of hyperplane arrangements is the relationship 
between combinatorics and topological properties. To be precise, one has to make the 
following distinction: for a given hyperplane arrangement $\cH\subset \PP^n$, one can 
study the topological type of the pair $(\PP^n,\cH)$ or the topological type of the 
complement $\PP^n\setminus \cH$. For the first concept we will use the term 
\emph{relative topology} of $\cH$, whereas for the second one we will simply say
\emph{topology} of $\cH$. It is clear that if two hyperplane arrangements have the same
relative topology, then they have the same topology, but the converse is not known.

In a well-known and very cited work~\cite{ry:98}, G.~Rybnikov found an example of two line
arrangements $L_1$ and $L_2$ in the complex projective plane $\PP^2$ having the same combinatorics but different topology. The most common way to prove that two topologies of 
line arrangements are different is to check that the fundamental groups of their complements
are not isomorphic.

Recently, the authors of this work have provided an example of two line arrangements with
different relative topologies (see~\cite{accm}). The contribution of~\cite{accm} is that 
it refers to real arrangements, that is, arrangements that admit real equations for each 
line (note that Rybnikov's example cannot have real equations).

The proof proposed by Rybnikov has two steps.
Let $G_i:=\pi_1(\PP^2\setminus\bigcup L_i)$, $i=1,2$.
\begin{enumerate}[(R1)]
\item \label{R1}
Recall that the homology of the complement of a hyperplane arrangement depends only on
combinatorics. This way, one can identify the abelianization of $G_1$ and $G_2$ with an 
Abelian group $H$ combinatorially determined. Rybnikov proves that no isomorphisms exist
between $G_1$ and $G_2$ that induce the identity on $H$. In particular, this result proves 
that both arrangements have different relative topologies. The reason can be outlined
as follows: any automorphism of the combinatorics of Rybnikov's arrangement can be 
obtained from a diffeomorphism of $\PP^2$, and hence it produces an automorphism of
fundamental groups. Since any homeomorphism of pairs $(\PP^2,\bigcup L_i)$ induces an
automorphism of the combinatorics of $\bigcup L_i$, after composition one can assume that 
any homeomorphism of pairs induces the identity on~$H$.

The strategy rests on the study of the Lower Central Series (LCS). Since $L_1$ and $L_2$ 
are constructed using the MacLane arrangement $L_\omega$ (see Example~\ref{mclane-rlz}), 
it is enough to study, by some combinatorial arguments, the LCS of $L_\omega$ with
an extra structure (referred to as an ordered arrangement). Although this part is explained
in~\cite[Section 3]{ry:98}, computations are hard to verify.

\item \label{R2}
The second step is essentially combinatorial. The main point is to truncate the LCS of
$G_i$ such that the quotient $K$ depends only on the combinatorics. Rybnikov proposes to
prove that an automorphism of $K$ induces the identity on $H$ (up to sign and automorphisms 
of the combinatorics). This proof is only outlined in~\cite[Proposition 4.2]{ry:98}. 
It is worth pointing out that such a result cannot be expected for any arrangement.
Also~\cite[Proposition 4.3]{ry:98} needs some explanation of its own. The main difference
between relative topology and topology of the complement in terms of isomorphisms of the
fundamental group is that homeomorphisms of pairs induce isomorphisms that send meridians 
to meridians, whereas homeomorphisms of the complement can induce any kind of isomorphism,
and even if we know that the isomorphism induces the identity on homology, this is not
enough to claim that meridians are sent to meridians.
\end{enumerate}

The aim of our work is to follow the idea behind Rybnikov's work and, using slightly 
different techniques, provide detailed proofs of his result. This is part of a more general
project by the authors that aims to better understand the relationship between topology 
and combinatorics of line arrangements.

The following is a more detailed description of the layout of this paper. In 
section~\ref{sec-defs}, the more relevant definitions are set, as well as a description of
Rybnikov's and MacLane's combinatorics. Sections~\ref{sec-alex}~and~\ref{sec-method} provide 
a proof of Step~(R\ref{R1}). In order to do so, we propose a new approach related to Derived
Series, which is also useful in the study of Characteristic Varieties and the Alexander
Invariant. The Alexander Invariant of a group $G$, with a fixed isomorphism 
$G/G'\approx \ZZ^r$, is the quotient $G'/G''$ considered as a module over the ring
$\gL:=\ZZ[t_1^{\pm1},\dots,t_r^{\pm1}]$, which is the group algebra of $\ZZ^r$. By fixing 
an ideal $\fm\subset\gL$, the truncated module by $\gL/\fm^j$ is studied. The problem is
reduced to solving a system of linear equations.

Section~\ref{sec-combinat} is devoted to the study of combinatorial properties of a line
arrangement which ensure that any automorphism of the fundamental group of the complement 
\emph{essentially} induces the identity on homology (that is, the analogous 
of~\cite[Proposition 4.2]{ry:98}). This is an interesting question that can be applied to
general line arrangements. For the sake of simplicity, we only present our progress on line
arrangements with double and triple points. This provides a proof for the second
step~(R\ref{R2}).

\section*{Acknowledgments}
We want to thank G.Rybnikov for his contribution and his support.
His work has been a challenge for us for years.

\section{Settings and Definitions}
\label{sec-defs}
In this section, some standard facts about line combinatorics and ordered line
combinatorics will be described. Special attention will be given to MacLane and
Rybnikov's line combinatorics.

\begin{definition}
A \emph{combinatorial type} (or simply a \emph{(line) combinatorics}) is a couple
$\scc:=(\cL,\cP)$, where $\cL$ is a finite set and $\cP\subset \cP(\cL)$, satisfying that:
\begin{enumerate}
\item 
For all $P\in\cP$, $\#\cP\geq 2$;
\item
For any $\ell_1,\ell_2\in \cL$, $\ell_1\neq\ell_2$, $\exists!P\in \cP$ such that 
$\ell_1,\ell_2\in P$. 
\end{enumerate}
An {\it ordered combinatorial type} $\scc^{\text{ord}}$ 
is a combinatorial type where $\cL$ is an ordered set.
\end{definition}

\begin{notation}
\label{not-comb-type}
Given a combinatorial type $\scc$, the multiplicity $m_P$ of $P\in\cP$ is the number 
of elements $L\in\cL$ such that $P\in L$; note that $m_P\geq 2$.
The \emph{multiplicity of} a combinatorial type is the number
$1-\# \cL + \sum_{P\in\cP} (m_P-1)$.
\end{notation}

\begin{figure}[hb]
\includegraphics[scale=.5]{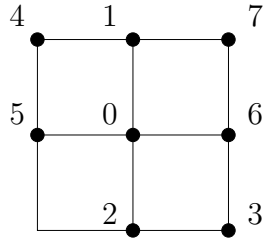}
\begin{picture}(0,0)
\put(-90,80){$4$}
\put(-90,43){$5$}
\put(-55,80){$1$}
\put(-55,43){$0$}
\put(-55,5){$2$}
\put(0,80){$7$}
\put(0,43){$6$}
\put(0,5){$3$}
\end{picture}
\caption{Ordered MacLane \emph{lines} in $\FF_3^2$}
\label{mclane}
\end{figure}

\begin{example}[MacLane's combinatorics]
\label{exam-maclane}
Let us consider the 2-dimensional vector space on  the field $\FF_3$ of three elements.
Such a plane contains 9 points and $12$ lines, $4$ of which pass through the origin. 
Consider $\cL=\FF_3^2\setminus \{(0,0)\}$ and $\cP$, the set of lines in $\FF_3^2$ 
(as a subset of $\cP(\cL)$).
This provides a combinatorial type structure $\mlc$ that we will refer 
to as \emph{MacLane's combinatorial type}. Figure~\ref{mclane} represents an ordered 
MacLane's combinatorial type.
\end{example}

\begin{definition}
Let $\scc:=(\cL,\cP)$ be a combinatorial type. We say a complex line arrangement 
$\cH:=\ell_0\cup \ell_1\cup ... \cup \ell_r\subset \PP^2$ is a 
\emph{realization} of $\scc$ if and only if there are bijections 
$\psi_1:\cL\to \{\ell_0,\ell_1,...,\ell_r\}$ and
$\psi_2:\cP\to \Sing(\cH)$ such that $\forall \ell\in \cH, P\in \cP$, one has
$P\in\ell \Leftrightarrow \psi_1(\ell)\in\psi_2(P)$. If $\scc^{\text{ord}}$ is an 
ordered combinatorial type and the irreducible components of $\cH$ are also ordered, 
we say $\cH$ is an ordered realization if $\psi_1$ respects orders.
\end{definition}

\begin{notation} 
The space of all complex realizations of a line combinatorics $\scc$ is denoted 
by $\Sigma(\scc)$. This is a quasiprojective subvariety of $\PP^{\frac{r(r+3)}{2}}$, where
$r:=\#\scc$. If $\scc^{\text{ord}}$ is ordered, we denote by  
$\Sigma^{\text{ord}}(\scc)\subset(\check\PP^2)^r$ 
the space of all ordered complex realizations of $\scc^{\text{ord}}$.
\end{notation}

There is a natural action of $\pgl(3;\CC)$ on such spaces. 
This justifies the following definition.

\begin{definition} 
The \emph{moduli space} of a combinatorics $\scc$
is the quotient $\scm(\scc):=\Sigma(\scc)/\pgl(3;\CC)$.
The \emph{ordered moduli space} $\scm^{\text{ord}}(\scc)$ of an 
ordered combinatorics $\scc^{\text{ord}}$ is defined accordingly.
\end{definition}

\begin{example} 
\label{mclane-rlz}Let us consider the MacLane line combinatorics $\mlc$. It is well known that 
such combinatorics has no real realization and that $\#\scm(\mlc)=1$, however
$\#\scm^{\text{ord}}(\mlc)=2$. The following are representatives for 
$\scm^{\text{ord}}(\mlc)$:

\begin{equation}
\label{eq-ryb-eqs}
\begin{gathered}
\ell_0 =  \{x=0\}\quad
\ell_1  =  \{y=0\}\quad
\ell_2  = \{x=y\}\quad
\ell_3  =  \{z=0\}\quad
\ell_4  =  \{x=z\}\\
\ell^{\pm}_5  =  \{z+\omega y=0\}\quad
\ell^{\pm}_6  =  \{z+\omega y=(\omega+1) x\}\quad
\ell^{\pm}_7  =  \{(\omega+1)y+z=x\}
\end{gathered}
\end{equation}

where $\omega=e^{2\pi i/3}$.

We will refer to such ordered realizations as 
$$L_\omega:=\{\ell_0,\ell_1,\ell_2,\ell_3,\ell_4,\ell_5^+,\ell_6^+,\ell_7^+\}$$ and 
$$L_{\bar \omega}:=\{\ell_0,\ell_1,\ell_2,\ell_3,\ell_4,\ell_5^-,\ell_6^-,\ell_7^-\}.$$
\end{example}

\begin{remark}
Given  a line combinatorics $\scc=(\cL,\cP)$, the automorphism group $\Aut(\scc)$
is the subgroup of the permutation group of $\cL$ preserving $\cP$.
Let us consider an ordered line combinatorics $\scc^{\text{ord}}$. It is easily seen that 
$\Aut(\scc^{\text{ord}})$ acts on both
$\Sigma^{\text{ord}}(\scc^{\text{ord}})$ and $\scm^{\text{ord}}(\scc^{\text{ord}})$. 
Note also that 
$\scm(\scc^{\text{ord}}) \cong \scm^{\text{ord}}(\scc^{\text{ord}})/\Aut(\scc^{\text{ord}})$.
\end{remark}

\begin{example} 
The action of $\Aut(\mlc)\cong\pgl(2,\FF_3)$ on the moduli spaces 
is as follows: matrices of determinant $+1$ (resp. $-1$) fix (resp. exchange) 
the two elements of $\scm^{\text{ord}}(\mlc)$.
Of course complex conjugation also acts on $\scm^{\text{ord}}(\mlc)$ exchanging the 
two elements. From the topological point of view one has that:
\begin{itemize}
\item There exists a homeomorphism 
$(\PP^2,\bigcup L_\gw)\to(\PP^2,\bigcup L_{\bar \gw})$
preserving orientations on both $\PP^2$ and the lines. Such a homeomorphism does not 
respect the ordering.

\item There exists a homeomorphism 
$(\PP^2,\bigcup L_\gw)\to(\PP^2,\bigcup L_{\bar \gw})$ preserving 
orientations on $\PP^2$, but not on the lines. Such a homeomorphism respects the ordering.
\end{itemize}

Also note that the subgroup of automorphisms that preserve the set
$L_0:=\{\ell_0,\ell_1,\ell_2\}$ is isomorphic to $\Sigma_3$, since the 
vectors $(1,0)$, $(1,1)$ and $(1,2)$ generate $\FF_3^2$. We will denote by $L_+$ 
and $L_-$ the sets of $5$ lines such that $L_\omega=L_0\cup L_+$ and 
$L_{\bar\omega}=L_0\cup L_-$. Since any transposition of $\{0,1,2\}$ in 
$\mlc$ produces a determinant $-1$ matrix in $\pgl(2,\FF_3)$, one concludes from the previous 
paragraph that any transposition of $\{0,1,2\}$ induces a homeomorphism
$(\PP^2,\bigcup L_\gw)\to(\PP^2,\bigcup L_{\bar \gw})$ that exchanges $L_\gw$ and 
$L_{\bar \gw}$ as representatives of elements of $\scm^{\text{ord}}(\mlc)$
and globally fixes $L_0$.
\end{example}

\begin{example}[Rybnikov's combinatorics]
\label{exam-ryb-comb}
Let $L_\omega$ and $L_{\bar \omega}$ be ordered MacLane realizations as above,
where $L_0:=\{\ell_0,\ell_1,\ell_2\}$. Let us consider a projective transformation 
$\rho_\omega$ (resp. $\rho_{\bar\omega}$) fixing the initial ordered set $L_0$ (that is, 
$\rho(\ell_i)=\ell_i$ $i=0,1,2$) and such that $\rho_\omega L_\omega$ 
(resp. $\rho_{\bar\omega} L_{\bar \omega}$) 
and $L_\omega$ intersect each other only in double points outside the three common 
lines. Note that $\rho_\omega,\rho_{\bar\omega}$ can be chosen with real coefficients.

Let us consider the following ordered arrangements of thirteen lines:
$R_{\alpha,\beta}=L_\alpha \cup \rho_{\gamma} L_\beta$, 
where
$\alpha,\beta\in \{\omega,\bar \omega\}$
and $\gamma=
\beta$ (resp $\bar\beta$) if $\alpha=\omega$
(resp. $\bar\omega$).
They produce the following combinatorics 
$\rybc:=(\cR,\cP)$ given by:
\begin{equation}
\label{eq-ryb-ct-2}
\begin{gathered}
\cR:=\{\ell_0,\ell_1,\ell_2,\ell_3,\ell_4,\ell_5,\ell_6,\ell_7,\ell_8,\ell_9,\ell_{10},
\ell_{11},\ell_{12}\}\\
\cP_2:=\left\{
\array{llll}
\{\ell_2,\ell_3\}, & \{\ell_0,\ell_7\}, & \{\ell_1,\ell_6\}, & \{\ell_4,\ell_5\}, \\
\{\ell_2,\ell_8\}, & \{\ell_0,\ell_{12}\}, & \{\ell_1,\ell_{11}\}, & \{\ell_9,\ell_{10}\}, \\
\{\ell_i,\ell_j\} & 3\leq i \leq 7, & 8\leq j \leq 12 \\
\endarray
\right\}\\
\cP_3:=\left\{
\array{lllll}
\{\ell_0,\ell_1,\ell_2\}, & \{\ell_3,\ell_6,\ell_7\}, & \{\ell_0,\ell_5,\ell_6\}, & 
\{\ell_1,\ell_4,\ell_7\}, & \{\ell_1,\ell_3,\ell_5\}, \\
\{\ell_2,\ell_4,\ell_6\}, & \{\ell_2,\ell_5,\ell_7\}, & \{\ell_0,\ell_3,\ell_4\}, &
\{\ell_8,\ell_{11},\ell_{12}\} & \{\ell_0,\ell_{10},\ell_{11}\}, \\
\{\ell_1,\ell_9,\ell_{12}\}, & \{\ell_1,\ell_8,\ell_{10}\}, & \{\ell_2,\ell_9,\ell_{11}\}, &
\{\ell_2,\ell_{10},\ell_{12}\}, & \{\ell_0,\ell_8,\ell_9\}
\endarray
\right\}\\
\cP:=\cP_2\cup \cP_3
\end{gathered}
\end{equation}
\end{example}

\begin{proposition} 
The following combinatorial properties hold:
\begin{enumerate}[\rm(1)]

\item The different arrangements $R_{\alpha,\beta}$ have the same combinatorial type
$\rybc$.

\item The set of lines $L_0$ has the following distinctive combinatorial property:
every line in $L_0$ contains exactly 5 triple points of the arrangement; the
remaining lines only contain 3 triple points. 

\item For the other $10$ lines we consider the equivalence relation generated by the 
relation of \emph{sharing a triple point}. There are two equivalence classes which correspond 
to $L_\varepsilon$ and $\rho L_{\varepsilon'}$, $\varepsilon,\varepsilon'=\pm$.

\end{enumerate}
\end{proposition}

By the previous remarks one can group the set $\cR$ together in three subsets.
One is associated with the set of lines $L_0$ (referred to as $\cR_0$), and
the other two are combinatorially indistinguishable sets ($\cR_1$ and $\cR_2$) such that 
$\cR_0 \cup \cR_1$ and $\cR_0 \cup \cR_2$ are MacLane's combinatorial types. Note 
that any automorphism of $\rybc$ must preserve $\cR_0$ and 
either preserve or exchange $\cR_1$ and $\cR_2$. Therefore,
$\Aut(\rybc)\cong\Sigma_3\times\ZZ/2\ZZ$. The following results are immediate
consequences of the aforementioned remarks.

\begin{proposition} 
The following are (or induce) homeomorphisms between the pairs
$(\PP^2,\bigcup R_{\bar \omega,\bar \omega})$ and 
$(\PP^2,\bigcup R_{\omega,\omega})$ 
(resp. $(\PP^2,\bigcup R_{\omega,\bar \omega})$ and $(\PP^2,\bigcup R_{\bar \omega,\omega})$) 
preserving the orientation of $\PP^2$:
\begin{enumerate}[\rm(a)] 

\item Complex conjugation, which reverses orientations of the lines.

\item A transposition in $\cR_0$, which preserves orientations of the lines.
\end{enumerate}
We will refer to $R_{\omega,\omega}$ and 
$R_{\bar \omega,\bar \omega}$ 
(resp. $R_{\omega,\bar \omega}$ and $R_{\bar \omega,\omega}$) 
as a type~$+$ (resp. type~$-$) arrangements.
\end{proposition}

\begin{proposition}
Any homeomorphism of pairs between a type~$+$ and a type~$-$ arrangement should
lead (maybe after composing with complex conjugation) to an orientation-preserving 
homeomorphism of pairs between a type~$+$ and a type~$-$ arrangement.

If such a homeomorphism existed, there should be an orientation-preserving 
homeomorphism of ordered MacLane arrangements of type $L_{\omega}$ and $L_{\bar \omega}$.
\end{proposition}

The purpose of the next section will be to prove that there is no 
orientation-preserving homeomorphism of ordered MacLane arrangements of type 
$L_{\omega}$ and $L_{\bar \omega}$.

\vspace*{14pt}
\section{The truncated Alexander Invariant}
\label{sec-alex}

Even though the Alexander Invariant can be developed for general projective
plane curves, we will concentrate on the case of line arrangements. 
Let $\bigcup L\subset \PP^2$ be a projective line arrangement where 
$L=\{\ell_0,\ell_1,...,\ell_r\}$. Let us denote its complement 
$X:=\PP^2 \setminus \bigcup L$ and $G$ its fundamental group.
The derived series associated with this group is recursively defined as follows: 
$G^{(0)}:=G$, $G^{(n)}:=(G^{(n-1)})'=[G^{(n-1)},G^{(n-1)}]$, $n\geq 1$, where $G'$ 
is the derived subgroup of $G$, i.e. the subgroup generated by 
$[a,b]:=a b a^{-1} b^{-1}$, $a,b \in G$. Note that the consecutive quotients are Abelian. 
This property also holds for the lower central series defined as 
$\gamma_1(G):=G$, $\gamma_n(G):=[\gamma_{n-1}(G),G]$, $n\geq 1$. 
It is clear that $G^{(0)}=\gamma_1(G)$ and $G^{(1)}=\gamma_2(G)$. 

Since $H_1(X)=G/G'$, one can consider the inclusion
$G'\inj G$ as representing the universal Abelian cover $\tilde X$ of $X$,
where $\pi_1(\tilde X)=G'$, and therefore $H_1(\tilde X)=G'/G''$.

The group of transformations $H_1(X)=G/G'=\ZZ^r$ of the cover acts on $G'$. 
This results in an action by conjugation on $G'/G''=H_1(\tilde X)$, $G''=G^{(2)}$:
$$
\array{cccccl}
G/G' \times G'/G'' & \rightarrow & G'/G'' &\\
(g,[a,b]) & \mapsto & g*[a,b] \ \Mod G''& =[g,[a,b]] + [a,b],
\endarray
$$
where $a*b:=a b a^{-1}$. This action is well defined since $g\in G'$ implies 
$g*[a,b]\equiv [a,b]$ mod $G''$. Additive notation will be used for the operation in $G'/G''$.

This action endows the Abelian group $G'/G''$ with a $G/G'$-module structure, that is, 
a module on the group ring $\Lambda:=\ZZ[G/G']$. If $x\in G$, then $t_x$ denotes its
class in $\Lambda$. For $i=1,\dots,r$, we choose $x_i\in G$ a meridian of $\ell_i$ in 
$G$; the class $t_i:=t_{x_i}\in\gL$ does not depend on the particular choice of the
meridian in $\ell_i$. Note that $t_1,\dots,t_r$ is a basis of $G/G'\cong\ZZ^r$ and
therefore one can identify
\begin{equation}
\label{eq-lambda}
\Lambda:=\ZZ[G/G']=\ZZ[t_1^{\pm 1},...,t_r^{\pm 1}].
\end{equation}
This module is denoted by $M_L$ and is referred to as the \emph{Alexander Invariant} of $L$. 
Since we are interested in oriented topological properties of $(\PP^2,\bigcup L)$,
the coordinates $t_1,\dots,t_r$ are well defined.

\begin{remark} 
The module structure of $M_L$ is in general complicated. 
One of its invariants is the zero set of the Fitting ideals of the 
complexified Alexander Invariant of $L$, that is, 
$M_L^\CC:=M_L \otimes \left( \gL\otimes \CC \right)$. This sequence of invariants is 
called the sequence of characteristic varieties of $L$
introduced by A.~Libgober \cite{li:01}. These are subvarieties of the 
torus $(\CC^*)^r$; in fact, irreducible components of characteristic varieties are translated subtori \cite{ara:97}.
\end{remark}

Our approach in studying the structure of the $\gL$-module $M_L$ is via the associated
graded module by the augmentation ideal $\fm:=(t_1-1,...,t_r-1).$ In order to do so,
and to be able to do calculations, we need some formul{\ae} on this module relating 
operations in $G'/G''$. For the sake of completeness, these formul{\ae} are listed below.
However, since they are straightforward consequences of the definitions, their proof will 
be omitted. The symbol ``$\eqmap{2}$'' means that the equality is considered in $G'/G''$:

\begin{properties}
\label{prop-alex-mod}
\nil
\begin{enumerate}
\item\label{prop-comm} $[x,p] \eqmap{2} (t_x-1)p$ $\ \forall p\in G'$, 
\item\label{prop-prod2} 
$\displaystyle[x_1\cdot ... \cdot x_n,y_1\cdot ... \cdot y_m] \eqmap{2} 
\sum_{i=1}^n \sum_{j=1}^m T_{ij} [x_i,y_j]$,
where
$\displaystyle
T_{ij}=\prod_{k=1}^{i-1} t_{x_k} \cdot 
\prod_{l=1}^{j-1} t_{y_l}$.
\item\label{prop-comm2}
$[p_1\cdot ...\cdot p_n,x] \eqmap{2} -(t_x-1)(p_1+...+p_n)$ $\ \forall p_i\in G'$, 
\item
$[p_xx,p_yy] \eqmap{2} [x,y]+(t_x-1)p_y-(t_y-1)p_x$ $\ \forall p_x,p_y\in G'$.
\item\label{prop-cnj}
$\displaystyle[x_1^{\alpha_1}\cdot\ldots\cdot x_n^{\alpha_n},y_1^{\beta_1}
\cdot\ldots\cdot y_m^{\beta_m}] \eqmap{2}
\sum_{i=1}^n \sum_{j=1}^m T_{ij} \Big( [x_i,y_j]+\delta(i,j) \Big),$
where 
$$
\delta(i,j)=(t_{y_j}-1)[\alpha_i,x_i]-(t_{x_i}-1)[\beta_j,y_j].
$$
\item\label{jacobi}\emph{Jacobi relations:}
$$J(x,y,z):=(t_x-1)[y,z]+(t_y-1)[z,x]+(t_z-1)[x,y] \eqmap{2} 0.$$
\end{enumerate}
\end{properties}

Let us recall a well-known result on presentations of fundamental groups of
line arrangements, an immediate consequence of the Zariski-Van Kampen method:

\begin{proposition}
\label{propos-pres}
The group $G$ admits a presentation of the form 
$\langle \bar x; \bar W \rangle$,
where $\bar x:=\{x_1,...,x_r\}$, and
$\bar W:=\{W_1(\bar x),...,W_m(\bar x)\}$ $m\geq 0$, such that $x_i$ 
is a meridian of $\ell_i$, and $W_i(\bar x)\in G'$ $\forall i=1,...,m$.
\end{proposition}

\begin{definition}
Any presentation $\langle \bar x; \bar W(\bar x) \rangle$
of $G$ as in Proposition~\ref{propos-pres} will be called a
\emph{Zariski presentation} of $G$. The free group $\langle \bar x \rangle$ 
will be referred to as the \emph{free group associated with} the given presentation.
\end{definition}

\begin{notation}
For technical reasons it is important to consider the Alexander Invariant corresponding 
to the free group associated with a given presentation. Such a module will be denoted by 
$\tilde M_L$. The following is a standard presentation of the modules $\tilde M_L$ 
and $\tilde M_L$ in terms of a
Zariski presentation of $G$.
\end{notation}

\begin{proposition}
\label{propos-pres-alex}
Let $\langle \bar x; \bar W \rangle$ be a Zariski presentation of $G$ and let
$\FF:=\langle \bar x \rangle$ be its associated free group, then the module $\tilde M_L$ 
admits a presentation $\tilde \Gamma/{\mathcal J}$, where 
$$\tilde \Gamma:=\bigoplus_{1\leq i<j \leq r} [x_i,x_j] \gL$$
and ${\mathcal J}$ is the submodule of $\tilde \Gamma$ generated by the Jacobi relations 
{\rm (Property~\ref{prop-alex-mod}(\ref{jacobi}))}
$$J(i,j,k):=(t_i-1)x_{jk}+(t_j-1)x_{ki}+(t_k-1)x_{ij}.$$
Moreover, the module $M_L$ can be obtained as a quotient of $\tilde M_L$ as
$\tilde \Gamma/({\mathcal J}+{\mathcal W})$, where ${\mathcal W}$ is the submodule of 
$\tilde\Gamma$ generated by the relations $\bar W$.
\end{proposition}

\begin{example}
\label{exam-dt}
As an example of how to obtain ${\mathcal W}$ note that, if the lines
$\ell_i$ and $\ell_j$ in $L$ intersect in a double point, then there is a relation in $\bar W$ of 
type $[x_i^{\ga_i},x_j^{\ga_j}]$, where $\ga_i,\ga_j\in G$. Using 
Property~\ref{prop-alex-mod}(\ref{prop-cnj}), this relation can be written in $\tilde M_L$
as $x_{i,j}+(t_j-1)[\ga_i,x_i]-(t_i-1)[\ga_j,x_j]\in \cW$. 

Analogously, if the lines $\ell_i$, $\ell_j$ and $\ell_k$ in $L$ intersect at a triple point,
one obtains relations in $G$ of type $[x_i^{\ga_i},x_j^{\ga_j}x_k^{\ga_k}]$, 
where $\ga_i,\ga_j,\ga_k\in G$, which can be rewritten in $\cW$ as
$$x_{i,j}+(t_j-1)[\ga_i,x_i]-(t_i-1)[\ga_j,x_j]+
t_j \big(x_{i,k}+(t_k-1)[\ga_i,x_i]-(t_i-1)[\ga_k,x_k]\big).$$
\end{example}

Let $\fm$ be the augmentation ideal in $\gL$ associated with the origin, that is, the 
kernel of homomorphism of $\gL$-modules, $\gve:\gL\to\ZZ$, $\gve(t_i):=1$, where
$\ZZ$ has the trivial module structure.

One can consider the filtration on $M_L$ associated with $\fm$, that is, 
$F^i M_L:= \fm^i M_L$. The associated graded module 
$\gr M_L:=\oplus_{i=0}^\infty \gr^i M_L$, where 
$\gr^i M_L:=F^i M_L/F^{i+1} M_L$ is a graded module over
$\gr_\fm \gL:=\oplus_{i=0}^\infty F^i \gL/F^{i+1} \gL$.

Consider the ring $\gL_j:=\gL/\fm^j$, tensoring $\gL$ by successive powers of the 
ideal $\fm$. This allows one to define truncations of the Alexander Invariant.

\begin{definition}
\label{def-trunc-alex-mod}
The $\gL_j$-module, $M_L^j:=M_L\otimes_\gL \gL_j$ will be called the 
\emph{$j$-th truncated Alexander Invariant of $L$}.
The induced filtration is finite and will be denoted in the same way.
\end{definition}

\begin{example}
\label{exam-dt-trunc}
From Example~\ref{exam-dt}, it is easy seen that the relations in $M_L^2$
coming from double and triple points are as follows.
\begin{enumerate}
\item\label{d} 
If $\ell_i$ and $\ell_j$ intersect at a double point one has:
\begin{equation}
x_{ij}+(t_j-1)[\ga_i,x_i]-(t_i-1)[\ga_j,x_j]=0,
\end{equation}
\begin{equation}
\label{exam-dt-trunc-dobles2} (t_k-1)x_{ij}=0.
\end{equation}
\item\label{t} 
If $\ell_i$, $\ell_j$ and $\ell_k$ intersect at a triple point one has:
\begin{equation}
x_{i,j}+t_j x_{i,k}+(t_j-1)[\ga_i,x_i]-(t_i-1)[\ga_j,x_j]+
(t_k-1)[\ga_i,x_i]-(t_i-1)[\ga_k,x_k]=0,
\end{equation}
\begin{equation}
\label{exam-dt-trunc-triples2}(t_m-1)x_{i,j}+(t_m-1)x_{i,k}=0.
\end{equation}
\end{enumerate}
\end{example}

For any $k\in \NN$, there is a natural morphism $\varphi_k:G'\surj M^k_L$. We will
sometimes refer to $\varphi_k(g)$ as $g \ (\Mod \fm^k)$ and equalities in $M^k_L$ will
be denoted by $p_1\equivmap{k}p_2$.

\begin{remark}
A Zariski presentation on $G$ induces a (set-theoretical) section 
$$
s:M_L\to G',
\quad
s\left( (t_1-1)^{k_1}...(t_r-1)^{k_r}x_{ij} 
\right):=[x_1,[x_i,x_j]]^{k_1}...[x_r,[x_i,x_j]]^{k_r},
$$
which, accordingly, induces a section of $\varphi_k$ on each $M^k_L$ denoted by $s_k$.
\end{remark}

\begin{remark}
From Property~\ref{prop-alex-mod}\eqref{prop-comm}, we deduce that the kernel of 
$\varphi_2$ equals $\gamma_4(G)$. 
Moreover $\ker(G'\to M_L^1)$ equals $\gamma_3(G)$.
\end{remark}

\begin{proposition}
\label{propos-modk}
Let $\psi(p_1,...,p_m)$ be a word on the letters
$\bar p:=\{p_1,...,p_m\}$. If $p_i,q_i\in G'$ and $p_i\equivmap{k}q_i$ $(i=1,...,m)$, then 
$[g,\psi(\bar p)]\equivmap{k+1}[g,\psi(\bar q)]$, $\forall g\in G$.
In particular, if $p\in M^k_L$ then $[g,p]$ is a well-defined
element of $M^{k+1}_L$; if $g=x_i$ this element
can be written $(t_i-1)p\in M^{k+1}_L$.
\end{proposition}

\begin{remark}
The ring $\gL_k$ is not local, but note that an element 
$\gl\in \gL_k$ is a unit if an only if $\gve(\gl)=\pm 1$. To see this 
note that $\gL_k=\ZZ \oplus \fm/\fm^k$ and the kernel of the evaluation
map $\gve:\ZZ \oplus \fm/\fm^k \to \ZZ$ is exactly $\fm/\fm^k$.
\end{remark}

\begin{notation}
\label{free-tilde} 
Note that everything in this section can also be 
reproduced by using the free group $\FF$ 
associated with a Zariski presentation of $G$ and will be denoted by adding a tilde. 
For instance, $\tilde M_L=\Lambda^{\binom{r}{2}}/{\mathcal J}$ 
is the Alexander Invariant associated with $\FF_G$, and $F^i \tilde M_L$ is the filtration 
associated with $\fm \subset \Lambda$.
\end{notation}

Note that the module $M_L^k$ is isomorphic (as an Abelian group) to 
the graduate 
$\gr^0 M_L^k\oplus ...\oplus \gr^{k-1} M_L^k$ by means of the  morphism
$$
p\cdot x_{ij}\mapsto p_0x_{ij} + p_1x_{ij} + ... + p_{k-1}x_{ij},
$$
where $p$ is a Laurent polynomial in $\bar t=\{t_1,...,t_r\}$ and 
$p_0+...+p_{k-1}$ is the $k$-truncated Taylor polynomial of $p$ around the origin
$\uno:=(1,...,1)$. Note that this isomorphism depends strongly on the Zariski 
presentation of $G$. More specifically, it depends on changes in the set of generators.

For instance, any automorphism of $G$ that sends $x_i$ to $x_i \alpha_i$, 
(with $\alpha_i\in G'$) induces both an automorphism of $M_L^k$ and a filtered 
automorphism of $\gr M_L^k$ which is the identity on $\gr^i M_L^k$, but 
not on $M_L^k$:
\begin{equation}
\label{eq-key}
[x_i,x_j]\mapsto [x_i {\alpha_i},x_j {\alpha_j}]\equivmap{k}
[x_i,x_j]+t_j(t_i-1)\alpha_j-t_i(t_j-1)\alpha_i.
\end{equation}

Nevertheless, the following result is true for these graduate groups.
\begin{proposition}
\label{prop-gr-comb}
The graduate groups $\gr^j M_L=\gr^j M_L^k$ (if $j<k$)
depend only on the combinatorics $\scc$ of $L$ and
will be denoted by $\gr^j M_\scc$.
\end{proposition}

The following result is an immediate consequence of Proposition~\ref{propos-modk} and 
it explains why $M_L^k$ is a more manageable object.

\begin{corollary}
\label{lineal}
The formula~{\rm(\ref{eq-key})} only depends on $\alpha_i\ (\Mod \fm^{k-1})$.
\end{corollary}

\section{Truncated Alexander Invariant and Homeomorphisms of Ordered Pairs}
\label{sec-method}

Let $L_1$ and $L_2$ be two ordered line arrangements
sharing the ordered combinatorics $\scc$. Consider two Zariski presentations 
$G_1=\langle \bar x;\bar W^1(\bar x) \rangle$ and
$G_2=\langle \bar x;\bar W^2(\bar x) \rangle$
of the fundamental groups of $X_{L_1}$ and $X_{L_2}$, where the subscripts of the 
generators $\bar x:=\{x_1,...,x_r\}$ respect the ordering of the irreducible 
components. 
The Abelian groups $G_1/G'_1$ and $G_2/G'_2$ can be canonically identified 
with $\gr^0 M_\scc$ so that 
$x_i\ (\Mod G'_1)\equiv x_i\ (\Mod G'_2)$. Hence $\gL:=\gL_{L_1}=\gL_{L_2}$
We will study the existence of isomorphisms $h:G_1\to G_2$ such that
$h_*:\gr^0 M_\scc\to \gr^0 M_\scc$ is the identity.

\begin{definition} 
Let $\FF_i$ be the free group associated with the Zariski presentation of $G_i$, $i=1,2$.
A morphism $\tilde h:\FF_1 \to \FF_2$ is called a \emph{homologically trivial morphism} if 
$\tilde h_*:\FF_1/\FF'_1\to\FF_2/\FF'_2$ satisfies $\tilde h_*(x_i)=x_i$. 
A morphism $h:G_1\to G_2$ 
is called a \emph{homologically trivial isomorphism} if it is induced by a homologically 
trivial morphism $\tilde h$, i.e., if $h_*:\gr^0 M_\scc\to \gr^0 M_\scc$
is the identity. Note that $\tilde h$ might not be unique.
\end{definition}

\begin{remarks}
\label{rem-htm}
\nil
\begin{enumerate}
\item\label{rem-htm-expl} 
In other words, a morphism $h:G_1\to G_2$ is homologically trivial if there exists
$(\alpha_1,...,\alpha_r)\in (G'_2)^r$ such that $h(x_i)=x_i \alpha_i$.

\item Any homologically trivial isomorphism $h$ induces a $\gL$-module morphism 
$h:M_1:=M_{L_1}\to M_{L_2}=:M_2$.

\item\label{rem-htm3}
Any homologically trivial isomorphism $h$ respects the filtrations $F$ and 
produces isomorphisms $\gr^i h:\gr^i M_1 \to\gr^i M_2$. By identifying 
$\gr^i M_1\equiv\gr^i M_\scc\equiv\gr^i M_2$, $\gr^i h$ is the \emph{identity}. 
\end{enumerate}
\end{remarks}

In order to state some properties of homologically trivial isomorphisms, we need to 
introduce some notation. Note that the homologically trivial morphism $\tilde h$ also 
induces morphisms on the Alexander Invariants of the associated free groups 
$\tilde M_i$ ($i=1,2$) and on their truncations $\tilde M^j_i$. 
Let us denote by $\tilde h^i:\tilde M_1^i\to \tilde M_2^i$ the induced homologically 
trivial morphisms of the truncated modules. A straightforward computation proves that
\begin{equation}
\label{prop-jacobi}
\tilde h (J(x_i,x_j,x_k))=J(x_i,x_j,x_k)\in\FF'_2/\FF''_2.
\end{equation}

Homologically trivial isomorphisms induce a particular kind of isomorphisms of the
$\gL$-modules $M_1,M_2$ which are worth studying.

\begin{remark}
\label{rem-gr-combinat}
A direct attempt to prove that two modules are homologically trivial isomorphic
is almost intractable. One would have to check if, for some choice 
$(\alpha_1,\dots,\alpha_r)\in (G'_2)^r\ \Mod G''_2$, such a $\gL$-module isomorphism exists.
The lack of linearity in this approach is the reason why we consider the truncated modules 
$M_1^k,M_2^k$.

Applying Corollary~\ref{lineal}, we are faced with simply solving a linear system as follows.
Let $h:G_1\to G_2$ be a homologically trivial morphism, then there exists
$(\alpha_1,\dots,\alpha_r)\in (G'_2)^r\ \Mod G''_2$ such that $h(x_i)=x_i\alpha_i$.
Therefore there exist $\Lambda_k$-morphisms $h^k: M^k_1\to M^k_2$ 
induced by $h$ for any $k\in \NN$. Note that 
$$h^2(x_{i,j}) = x_{i,j}+\sum_{u,v} \alpha^j_{u,v}x_{i,u,v}  - 
\sum_{u,v} \alpha^i_{u,v}x_{j,u,v}$$
where $x_{i,j}\equivmap{2}[x_i,x_j]$, $x_{i,u,v}\equivmap{2}(t_i-1) x_{u,v}$ and
$\alpha_w\equivmap{1}\sum_{u<v} \alpha^w_{u,v} x_{u,v}$, 
since
\begin{equation}
\label{eq-h2}
h^2(x_{i,j})\equivmap{2}[h(x_i),h(x_j)]\equivmap{2}[x_i\alpha_i,x_j\alpha_j]
\end{equation}
only depends on $\varphi_1(\alpha_i)$, the class of $\alpha_i\ \Mod \fm$, by 
Proposition~\ref{lineal}.
In order to prove that $h^2$ is well defined, one must solve
a linear system of equations on the variables $\alpha^i_{u,v}$ in the Abelian
group $M^k_2$. If an integer solution exists, one can repeat the procedure on $M_i^3$, 
obtaining again a linear system of equations in the Abelian group $M_2^3$, and so on.
In this work, we only need to consider $h^2$.
\end{remark}

Let us consider an ordered line arrangement $L$ with a fixed Zariski
presentation $G=\langle\bar x;\bar W\rangle$. Let us denote $\scc$ its
ordered combinatorics.

\begin{lemma}
\label{free}
The $\gL_1$-module $\gr^0 M_L={\gL_1^{\binom{r}{2}}}/{\mathcal W}$ is free of rank $g=\binom{r}{2} -v$, 
where $v$ is the multiplicity of the combinatorial type of $L$ 
(Notation~{\rm\ref{not-comb-type}}).
\end{lemma}

\begin{proof}
It is straightforward by taking suitable orderings on the relations.
\end{proof}

\begin{lemma}
\label{jac-homo}
Let $\widehat{\mathcal J}:=({\mathcal J}+{\mathcal W})/{\mathcal W}$. 
Then $\widehat{\mathcal J}\subset F^1\left(\gL_2^{\binom{r}{2}}/\cW\right)$.
\end{lemma}

\begin{proof}
A consequence of the fact that $J(i,j,k)\in F^1\gL_2^{\binom{r}{2}}$.
\end{proof}

\begin{remark}
\label{rem-gen-gr1}
Note that $F^1 M_L^2=\gr^1 M_L=\gr^1 M_\scc$ is an Abelian group with \emph{trivial} 
action of $\Lambda_2$. A system of generators is given by $(t_i-1)x_{j,k}$; the relations 
as an Abelian group are exactly those in $\mathcal J$ and the relations~\eqref{exam-dt-trunc-dobles2} and~\eqref{exam-dt-trunc-triples2}
in Example~\ref{exam-dt-trunc},
which depend only on the combinatorial type by~\eqref{eq-key}. 
\end{remark}

In our particular case we can effectively compute the 2-nd truncated Alexander Invariant.  
The following result is an easy computation.

\begin{lemma} 
\label{lem-mcl-rank}
For any MacLane arrangement $L$, the Abelian group $M_{L}^2$ is free of rank~$29$, 
and its subgroup $\gr^1 M_{L}^2=\gr^1 M_{\mlc}$ is free of rank~$21$.
\end{lemma}

\begin{theorem}
\label{thm-main}
There is no homologically trivial isomorphism between 
$G_\omega:=\pi_1(\PP^2\setminus \bigcup L_\omega)$ and 
$G_{\bar \omega}:=\pi_1(\PP^2\setminus \bigcup L_{\bar \omega})$.
\end{theorem}

\begin{proof}
Consider suitable Zariski presentations 
$G_\omega=\langle x_1,...,x_7;W^\omega_1(\bar x),...,W^\omega_{13}(\bar x) \rangle$ and 
$G_{\bar \omega}=\langle x_1,...,x_7;
W^{\bar \omega}_1(\bar x),...,W^{\bar \omega}_{13}(\bar x) \rangle$
of the ordered arrangements $L_\omega$ and $L_{\bar \omega}$ (for instance, we have used
the suitable Zariski presentations provided in~\cite{ry:98} and other presentations
obtained using the software in~\cite{car:xx}).
We identify the corresponding free groups $\FF_\omega$ and $\FF_{\bar\omega}$
with a free group $\FF_7$.
Recall that their combinatorial type has multiplicity $13$ and we suppose that the relations 
are ordered upon the following condition:
$$
W^\omega_i(\bar x)\equivmap{1}W^{\bar \omega}_i(\bar x),$$
(in particular 
$W^\omega_i(\bar x)-W^{\bar \omega}_i(\bar x)\in F^1 \tilde M_{L_{\bar\omega}}$).
Let us suppose that a homologically trivial homomorphism 
$h:G_\omega\to G_{\bar \omega}$ exists. Consider the corresponding 
elements $\alpha_1,\dots,\alpha_r\in\FF_7$ that induce such a morphism
(Remark~\ref{rem-htm}\eqref{rem-htm-expl}). Consider $M_{\omega}=M_{L_\omega}$ and 
$M_{\bar \omega}=M_{L_{\bar \omega}}$, the Alexander Invariants of 
${L_\omega}$ and $L_{\bar \omega}$. 
Let $\gL:=\ZZ[t_1^{\pm 1},...,t_7^{\pm 1}]$ be the ground ring of both Alexander Invariants, 
where $t_i\eqmap{1} x_i$ as usual. 
This mapping induces a 
$\gL_2$-isomorphism $h^2:M_{\omega}^2\to M_{\bar \omega}^2$. By Corollary~\ref{lineal}, 
$h^2$ only depends on the class $\alpha_i\mod \fm$. One has 
$$\alpha_k\equivmap{1}\sum_{j=1}^7 p^k_{ij} x_{ij},\quad p^k_{ij}\in\ZZ.$$ 
By \eqref{prop-jacobi} Jacobi relations play no role here.

Let us fix $i=1,\dots,13$. Since 
$W^\omega_i(\bar x)\in \tilde M_{\omega}^2$ vanishes in $M_{\omega}^2$, one deduces that
$\tilde h^2(W^\omega_i(\bar x))\in \tilde M_{L_{\bar\omega}}^2$ should vanish in
$M_{\bar\omega}^2$. Equivalently
$\tilde h^2(W^\omega_i(\bar x))-W^{\bar \omega}_i(\bar x) \in F^1 \tilde M_{\bar\omega}^2$
should also vanish in $F^1 M_{\bar\omega}^2=\gr^1 M_{\mlc}$. The vanishing of these terms, 
considered in the free group 
$\gr^1 M_{\mlc}$, produces a system of linear equations in the variables $p^k_{ij}$
(actually, even though there are 147 variables, only 126 appear in the equations).

Solving a system with 137 equations and 126 variables is not an easy task, but any computer
will help. Using Maple8, it takes 85 seconds of CPU time running on an Athlon at 1.4MHz and
256Kb RAM Memory to obtain the linear set of solutions. It is an affine variety of dimension
98 of the form $(\gl_1,...,\gl_{98},\gk_1,...,\gk_{28})$ where 
$$\gk_i=q_i+\sum_{j=1}^{98} \gve_j^i \gl_j,$$
$\gve_j^i\in \{0,\pm 1\}$, and $q_i\in \QQ$. Since 
$\{q_i\mid i=1,\dots,28\}=
\{0,\pm 1,\pm \frac{1}{3},\pm \frac{2}{3},\pm \frac{4}{3},\pm \frac{5}{3}\}$ 
one concludes that there is no integer solution\footnote{The software is written for
Maple8 and can be visited at the following public site
\texttt{http://riemann.unizar.es/geotop/pub/.}}.
\end{proof}

\begin{corollary}
There is no orientation preserving homeomorphism between the pairs of ordered arrangements 
$(\PP^2,\bigcup L_\omega)$ and $(\PP^2,\bigcup L_{\bar \omega})$.
\end{corollary}

For Rybnikov's arrangements, one obtains similar results.

\begin{lemma} 
For any Rybnikov's arrangement $R$, the Abelian group $M_{R}^2$ is free of rank~$55$,
and its subgroup $\gr^1 M^2_{R}=\gr^1 M_{\rybc}$ is free of rank~$40$.
\end{lemma}

\begin{theorem}
\label{thm-main-ryb}
There is no homologically trivial isomorphism between
$G_{+}:=\pi_1(\PP^2\setminus \bigcup R_{\omega,\omega})$ and 
$G_{-}:=\pi_1(\PP^2\setminus \bigcup R_{\bar \omega,\omega})$.
\end{theorem}

\begin{proof}
One way to prove this statement is to follow the computational strategy proposed for MacLane arrangements. 
First one needs Zariski presentations of $G_{\pm}$. This was done by
means of the software in~\cite{car:xx}. In this case the linear system
obtained consists of 531 equations and 420 variables (again, out of the 792
variables $p^k_{ij}$, only 420 appear in the equations) and it took the same
processor a total of 23,853 seconds of CPU time to compute the solutions. The space of solutions has dimension 252, that is, it can be written as  $(\gl_1,...,\gl_{252},\gk_1,...,\gk_{168})$, where 
$$\gk_i=q_i+\sum_{j=1}^{168} \gve_j^i \gl_j,$$
$\gve_j^i\in \{0,\pm 1\}$, and $q_i\in \QQ$.
Since $\{q_i\mid i=1,\dots,168\}=
\{0,\pm 1,\pm \frac{1}{3},\pm \frac{2}{3},\pm \frac{4}{3},\pm \frac{5}{3}\}$, one again
concludes that there is no integer solution.

Another proof that doesn't depend as strongly on computations can be obtained from 
Theorem~\ref{thm-main} as follows. Let us assume that a homologically trivial isomorphism 
exists between $G_{+}$ and $G_{-}$. Such an isomorphism induces an $\gL_2$-isomorphism
between $M^2_+$ and $M^2_-$. Let $\widehat \gL_2:=\gL_2/\fm'$, where 
$\fm'$ is the ideal generated by $(t_8-1),...,(t_{12}-1)$, and let
$\widehat M^2_\pm$ denote $M^2_\pm \otimes \widehat \gL_2$. Note that $M^2_{\gw}$ (resp. $M^2_{\bar \gw}$)
can be considered as the $\widehat \gL_2$-module obtained from the inclusion of 
the complements $\PP^2\setminus \bigcup R_{\omega,\omega}\inj \PP^2\setminus\bigcup L_{\gw}$
(resp. $\PP^2\setminus \bigcup R_{\bar \omega,\omega}\inj \PP^2\setminus\bigcup L_{\bar \gw}$).
Moreover, these inclusions define epimorphisms of
$\widehat \gL_2$-modules $\pi:\widehat M^2_+\twoheadrightarrow M^2_{\gw}$
and $\bar\pi:\widehat M^2_+\twoheadrightarrow M^2_{\bar\gw}$.
Proving the existence of a homologically trivial isomorphism
$\tilde h^2$ that matches in the commutative diagram~\eqref{eq-comm-omega},
and using Theorem~\ref{thm-main} one obtains a contradiction.
\begin{equation}
\label{eq-comm-omega}
\begin{matrix}
\widehat M^2_+ & \longrightmap{\hat h^2} & \widehat M^2_-\\
\pi\downarrow&&\downarrow\bar\pi
\\
M^2_{\gw}& \surjmap{\tilde h^2} & M^2_{\bar \gw}
\end{matrix}
\end{equation}
Consider $S_+$ the $\widehat \gL_2$-submodule of $\widehat M^2_+$ generated by the elements $x_{i,j}$, $i,j\in \{1,...,7\}$
and consider the commutative diagram~\eqref{eq-comm-somega}.
Since $\pi(x_{i,j})=x_{i,j}$ and 
$\bar h^2(x_{i,j})\equiv x_{i,j}\mod \gr^1 M^2_{\bar\omega}$.
\begin{equation}
\label{eq-comm-somega}
\begin{matrix}
\widehat M^2_+ & \longrightmap{\hat h^2} & \widehat M^2_-\\
\pi\ninj&&\downarrow\bar\pi
\\
S_+& \surjmap{\bar h^2} & M^2_{\bar \gw}\\
\pi_{|}\ssurj&&\\
M^2_{\gw}
\end{matrix}
\end{equation}
%
Since $M^2_{\gw}$ and $M^2_{\bar\gw}$ are free Abelian of the same rank,
\eqref{eq-comm-omega} can be obtained from \eqref{eq-comm-somega},
by proving that $\pi_|$ is an isomorphism,
which is the statement of Lemma~\ref{lema-comm-diag}.
\end{proof}

\begin{lemma}
\label{lema-comm-diag}
The epimorphism $\pi_{|}$ in~\eqref{eq-comm-somega} is injective
\end{lemma}

\begin{proof}
We break the proof in several steps.
\begin{enumerate}[(O1)]
\smallbreak\item\label{ceros}
$(t_k-1)x_{i,j}=0$ in $\widehat M^2_{+}$ if $\{i,j\}\cap \{8,...,12\}\neq \emptyset$.

This can be proved case by case (all the equalities are considered in $\widehat M^2_{+}$):
\makeatletter\renewcommand{\p@enumii}{}
\begin{enumerate}
\item\label{cerosa} If $i\in \{3,...,7\}$ and $j\in \{8,...,12\}$ 
(or vice versa, since $x_{i,j}=-x_{j,i}$): this is a consequence of
Example~\ref{exam-dt-trunc}\eqref{d} since the lines $\ell_i$ and $\ell_j$ intersect
transversally.
\item\label{cerosb} If $i,j\in \{8,...,12\}$: this is a consequence of (\ref{cerosa}) and the 
Jacobi relations (Property~\ref{prop-alex-mod}(\ref{jacobi})).
\item If $i\in \{1,2\}, j\in \{8,...,12\}$: using the Jacobi relations
(Property~\ref{prop-alex-mod}(\ref{jacobi})) and (\ref{cerosb})
it is enough to check that $(t_k-1)x_{i,j}=0$, 
$i,k\in \{1,2\}$, $j\in \{8,...,12\}$. If $\ell_i$ and $\ell_j$ intersect at a double point,
then (\ref{cerosa}) proves the result. Otherwise, there exists a line $\ell_m$ 
($m\in \{7,...,12\}$) such that $\ell_i$, $\ell_j$ and $\ell_m$ intersect at a triple point.
By Example~\ref{exam-dt-trunc}(\ref{t}) one has $(t_k-1)x_{i,j}+(t_k-1)x_{m,j}=0$, 
but $(t_k-1)x_{m,j}=0$ by (\ref{cerosb}), thus we are done.
\end{enumerate}

\smallbreak\item $\gr^1(\widehat M^2_+)\subset S_+$.
It is a direct consequence of (O\ref{ceros}).

\smallbreak\item\label{otres}
$\displaystyle\gr^0(\widehat M^2_+)=\frac{S_+}{\gr^1 \widehat M^2_+} \oplus 
\frac{\ker\pi+\gr^1 \widehat M^2_+}{\gr^1 \widehat M^2_+}$, i.e.,
$\dfrac{S_+}{\gr^1 \widehat M^2_+}\cong\gr^0 M^2_\gw$.

Since $\ker\pi$ is generated by $\langle x_{i,j} \rangle_{1\leq i < j\leq 12,j>7}$, it is clear that $\gr^0(\widehat M^2_+)$ decomposes
in the required sum. It remains to prove that it is a direct sum.

One can consider $\gr^0(\widehat M^2_+)$ as a quotient 
$\dfrac{\langle x_{i,j} \rangle_{1\leq i < j\leq 12}}{\cW}$, hence it is enough to check that there 
is a system of generators $r_1,...,r_n$ of $\cW$ such that:
\begin{equation}\label{condition}
\tag{*}
\text{either }\quad  
r_i\in \frac{\ker\pi+\gr^1 \widehat M^2_+}{\gr^1 \widehat M^2_+},\quad
\text{or }\quad r_i\in \dfrac{S_+}{\gr^1 \widehat M^2_+}.
\end{equation}
 Note that a system of relators can be obtained
combinatorially as $x_{i,j}=0$ (if $\{\ell_i,\ell_j\}$ is a double point) or
$x_{i,j}+x_{i,k}=0$ (if $\{\ell_i,\ell_j,\ell_k\}$ is a triple point). 
Relations coming from double points satisfy \eqref{condition}. For the
triple point relations note that any triple point
$\{\ell_i,\ell_j,\ell_k\}$ such that $\{i,j\}\subset \{8,...,12\}$, verifies that 
$k\in \{8,...,12\}$; therefore, condition~\eqref{condition}
is also satisfied.

%

\smallbreak\item $\gr^1 (\widehat M^2_+)\cong\gr^1 M^2_\gw$.

By (O\ref{ceros}), the Abelian group $\widehat M^2_+$ is generated by
$x_{i,j}$, $i,j\in\{1,\dots,12\}$, and
$(t_k-1) x_{i,j}$, $i,j,k\in\{1,\dots,7\}$; the relators are obtained
from the singular points (see Example~\ref{exam-dt-trunc}) and the Jacobi relations $\mathcal J$. 

By (O\ref{ceros}) and Remark~\ref{rem-gen-gr1},
we find that $\gr^1 (\widehat M^2_+)$ is generated by the elements
$(t_k-1) x_{i,j}$, $i,j,k\in\{1,\dots,7\}$ and the relations
are exactly those in $\mathcal J$ and the relations~\eqref{exam-dt-trunc-dobles2} and~\eqref{exam-dt-trunc-triples2}
in Example~\ref{exam-dt-trunc}. The arguments used in (O\ref{otres})
also show that only double and triple points in $\mlc$ provide non-trivial
relations and thus one obtains the same system of generators and 
relations of $\gr^1 M^2_\gw$.
\end{enumerate}
\end{proof}

\begin{remark}
The proof of Lemma~\ref{lema-comm-diag} is combinatorial
and depends strongly on the properties of $\rybc$. This lemma
corresponds to a key statement of the proof of
\cite[Lemma 4.3]{ry:98} which is worth mentioning.
\end{remark}

\section{Homologically Trivial Fundamental Groups}
\label{sec-combinat}

This last section will be devoted to proving that the fundamental groups of $R_{\gw,\gw}$
and $R_{\gw,\bar \gw}$ are not isomorphic.

\begin{remark}
\label{rem-m-combinat}
Associated with a combinatorial type $\scc:=(\cL,\cP)$, there is a family of groups, 
where
$$H_\scc:=\frac{\bigoplus_{\ell\in \cL} \langle x_{\ell} \rangle \ZZ}
{\langle \sum_{\ell\in \cL} x_\ell \rangle \ZZ}$$
and $\gr^iM_\scc$ is given by generators and relations as a quotient of 
$H_\scc^{\otimes i}$ as described in Remark~\ref{rem-gr-combinat}. 
Note that, if $\scc$ has a realization $L$, then one has identifications
$H_\scc \equiv H_1(\PP^2\setminus \bigcup L;\ZZ)$ and 
$\gr^iM_\scc\equiv\gr^i M_L$.
\end{remark}

\begin{notation}
There is a natural injective map $\Gamma(\scc)\hookrightarrow \Aut(H_\scc)$ given by 
the permutation of the generators of $H_\scc$;
we identify $\Gamma(\scc)$ with its image in $\Aut(H_\scc)$.
Another subgroup of 
$\Aut (H_\scc)$, denoted by $\Aut^1(H_\scc)$, is defined as those automorphisms of
$H_\scc$ that induce an automorphism of $\gr^1 M_\scc$.
It is easily seen that $\{\pm 1_{H_\scc}\}\times\Aut(H_\scc)\subset\Aut^1(H_\scc)$.
\end{notation}

\begin{definition}
A line combinatorics $\scc:=(\cL,\cP)$ is called \emph{homologically rigid} if
$\Aut^1(H_\scc)=\{\pm 1\} \times \Gamma(\scc)$.
\end{definition}

The first goal of this section is to prove that Rybnikov's combinatorial type 
$\rybc:=(\cR,\cP)$ (described in~(\ref{eq-ryb-ct-2})) is homologically rigid;
we will follow the ordering~(\ref{eq-ryb-eqs}).
In order to do so, 
we are going to study $\Aut^1(H_\scc)$ for an ordered combinatorics $\scc=(\cL,\cP)$
having at most triple points. We will denote
$\cP_j:=\{P\in\cP\mid\#P=j\}$, $j=2,3$, and
$\cL:=\{\ell_0,\ell_1,\dots,\ell_r\}$.
Let us first describe the groups $H_\scc$ and $\gr^1M_\scc$:
$$H_\scc:=\frac{\langle x_0 \rangle \ZZ \oplus \langle x_1 \rangle \ZZ \oplus\dots
\oplus \langle x_{r} \rangle \ZZ}{\langle x_0+x_1+...+x_{r}\rangle \ZZ}, \quad
\gr^1M_\scc:=\frac{\bigwedge^2 H_\scc}{R_2 \oplus R_3}$$
where $R_2$ is the subgroup generated by $x_{i,j}$ 
($\{\ell_i,\ell_j\}\in \cP_2$) and $R_3$ is the subgroup generated by
$x_{i,k}+x_{j,k}$ and $x_{i,j}+x_{i,k}$ ($\{\ell_i,\ell_j,\ell_k\}\in \cP_3$).

Any isomorphism $\psi:H_\scc\to H_\scc$ induces a map 
$\wedge^2 \psi:\bigwedge^2 H_\scc \to \bigwedge^2 H_\scc$. 
Let us represent $\psi:H_\scc\to H_\scc$ by means of a matrix
$A^\psi:=(a^j_i)\in \Mat(r+1,\ZZ)$ such that $\psi(x_i):=\sum_{j=0}^{r} a_i^j x_j$
(note that such a matrix is not uniquely determined: each column is only well defined
modulo the vector $\uno_{r+1}:=(1,\dots,1)$). 
The conditions required for this map to define a morphism on the quotient $\gr^1M_\scc$ are 
called \emph{admissibility conditions} and can be expressed as follows:

\begin{equation}
\label{determinantes1}
\array{lll}
\left\vert
\begin{array}{ccc}
a_u^i &a_v^i &1\\
a_u^j &a_v^j &1\\
a_u^k &a_v^k &1 
\end{array}
\right\vert=0,
& \text{ if } \{\ell_i,\ell_j,\ell_k\}\in\cP_3,
\{\ell_u,\ell_v\}\in \cP_2 &\\
\\
\left\vert
\begin{array}{ccc}
a_\bullet^i &a_u^i+a_v^i+a_w^i &1\\
a_\bullet^j &a_u^j+a_v^j+a_w^j &1\\
a_\bullet^k& a_u^k+a_v^k+a_w^k &1
\end{array}
\right\vert=0
& \text{ if } \{\ell_i,\ell_j,\ell_k\},
\{\ell_u,\ell_v,\ell_w\}\in \cP_3&(\bullet=u,v,w)\\
\endarray
\end{equation}
(also note that such conditions are invariant on the coefficient vectors 
$(a_\bullet^0,a_\bullet^1,...,a_\bullet^{12})$ modulo $\uno_{13}$).
We summarize these facts.

\begin{proposition}
\label{prop-adm-cond}
Any morphism $\psi:H_\scc\to H_\scc$ whose associated matrix $A^\psi$ satisfies the 
admissibility conditions \eqref{determinantes1} produces a well-defined morphism 
$\wedge^2 \psi:M^1_\scc\to M^1_\scc$.
\end{proposition}

We are going to express the admissibility conditions \eqref{determinantes1}
in a more useful way. Let $\psi\in\Aut^1(H_\scc)$ and
let $A^\psi$ be a matrix representing $\psi$.
Fix $P\in \cP_3$ and consider the submatrix
$A^\psi_P\in \Mat(3\times 12,\ZZ)$ of $A^\psi$ which contains
the rows associated with $P$. Let $\Sigma_k:=\ZZ^{k+1}/\uno_{k+1}$, $k\in\NN$. 
We denote by $v_0(P),v_1(P)\dots,v_r(P)\in\Sigma_2$,
the column vectors $\!\!\!\pmod{\uno_3}$ of $A^\psi_P$.

\begin{lemma}
\label{lem-adm-cond}
\nil
\begin{enumerate}[\rm(1)]
\item\label{lem-adm-cond1} The vectors $v_0(P),v_1(P)\dots,v_r(P)$ span $\Sigma_2$.
\item \label{lem-adm-cond2}$\sum_{j=0}^r v_j(P)=0\in\Sigma_2$.
\item \label{lem-adm-cond3}For any $Q\in\cP$ and for any $\ell_u\in Q$,
the vectors $v_u(P)$ and $\sum_{\ell_i\in Q} v_i(P)$
are linearly dependent (i.e., span a sublattice of $\Sigma_2$ of rank less than two).
In particular, if $\sum_{\ell_i\in Q} v_i(P)\neq 0$, then
$\{v_i(P)\mid\ell_i\in Q\}$ spans a rank-one sublattice of $\Sigma_2$.
\item\label{lem-adm-cond4} There exists $Q\in\cP_3$ such that
$\{v_i(P)\mid\ell_i\in Q\}$ spans a rank-two sublattice of $\Sigma_2$
and $\sum_{\ell_i\in Q} v_i(P)=0$.
\end{enumerate}
\end{lemma}

\begin{proof}
\nil
\begin{itemize}
\item[\eqref{lem-adm-cond1}] $\psi$ is an automorphism.
\item[\eqref{lem-adm-cond2}] The sum of the columns of $A^\psi$
is a multiple of $\uno_{r+1}$.
\item[\eqref{lem-adm-cond3}] It
is an immediate consequence of the admissibility conditions \eqref{determinantes1}.
\item[\eqref{lem-adm-cond4}] If no such $Q$ exists, then all
the vectors $V_i(P)$ are linearly dependent, which contradicts
\eqref{lem-adm-cond1}. The last part is a consequence
of \eqref{lem-adm-cond3}.
\end{itemize}
\end{proof}

\begin{definition}
\label{def-subcomb}
Let $\scc:=(\cL,\cP)$ be a combinatorics; we say $\scc':=(\cL',\cP')$ is a 
\emph{subcombinatorics of $\scc$} if $\cL'\subset \cL$ and 
$\cP':=\{P\cap \cL\mid P\in \cP,\#(P\cap\cL)\geq 2\}$. 
\end{definition}

We define a subcombinatorics $\Adm_\psi(P)\subset \scc$ as follows:
$$
\cL(\Adm_\psi(P)):=\{\ell_i\in \cR \mid v_i(P)\neq 0\}.
$$
Note that, 
\begin{equation}
\label{eq-unos}
\ell_i\notin \cL(\Adm_\psi(P))\Longleftrightarrow 
\text{ the }i^{\text{th}}\text{ column of }A^\psi_P
\text{ is a multiple of }\uno_3. 
\end{equation}
This motivates the following definition.

\begin{definition}
\label{def-adm}
A line combinatorics $\scc:=(\cL,\cP)$ with only double and triple points
is called \emph{3-admissible} if it is possible to assign a non-zero vector 
$v_i\in \ZZ^2$ to each $\ell_i\in \cL$ such that:
\begin{enumerate}
\item\label{def-adm-1} 
There exists $P\in \cP_3$, such
that $\{v_j\mid \ell_j \in P\}$ spans a rank-two sublattice.
\item\label{def-adm-1b} 
For every $P\in\cP$ and for every $\ell_i\in P$,
$v_i$ and $\sum_{\ell_j \in P} v_j$ are linearly dependent.
\item\label{def-adm-3}
$\sum_{\ell_i \in \cL} v_i=(0,0)$.
\end{enumerate}
\end{definition}

\begin{remarks}
\label{ref-def-adm}
The conditions of Definition~\ref{def-adm}
can be made more precise.
\begin{enumerate}
\item\label{ref-def-adm2} If $P=\{\ell_i,\ell_j\}\in\cP_2$, then
$v_i$ and $v_j$ are proportional, in notation, $v_i\Vert v_j$.
\item\label{ref-def-adm3} If $P\in\cP_3$ verifies condition \eqref{def-adm-1}
then $\sum_{\ell_i \in P} v_i=(0,0)$.
\end{enumerate}
\end{remarks}

\begin{examples}
\label{exam-adm}
\nil
\begin{enumerate}
\item\label{exam-adm-adm} With the above notation, $\Adm_\psi(P)$
is 3-admissible by Lemma~\ref{lem-adm-cond}.
\item\label{exam-adm-triple}
The combinatorics $\cM_3$ of a triple point 
(that is, 
$\cL_{\cM_3}:=\{0,1,2\}$, $\cP_{\cM_3}:=\{\{0,1,2\}\}$)
is 3-admissible, simply using $v_0:=(1,0)$, $v_1:=(0,1)$, $v_2:=(-1,-1)$.
\item\label{exam-adm-general}
Let $\scc:=(\cL,\cP)$ be a combinatorics such that 
\begin{itemize}
\item $\cL:=\cL_0 \coprod \cL_1\coprod \cL_2$;
\item $\cL_1$ and
$\cL_2$ define non-empty subcombinatorics in general position w.r.t. $\cL_0$ 
(that is, $\ell_i\in \cL_1$ and $\ell_j\in \cL_2$ implies that $\{\ell_i,\ell_j\}\in \cP$);
\item at most one line of $\cL_0$ \emph{intersects} $\cL_1\cup\cL_2$
in a non-multiple point of $\cL_0$.
\end{itemize}
Then $\scc$ is not 3-admissible. It is enough to see that if 
$\{v_i\}_{\ell_i\in \cL}$ is a set of non-zero vectors satisfying
conditions~\eqref{def-adm-1b}~and~\eqref{def-adm-3}
of Definition~\ref{def-adm}, one has that $v_i\Vert v_j$. 

\begin{itemize}
\item $\ell_i\in \cL_1$ and $\ell_2\in \cL_2$;
since $\{\ell_i,\ell_j\}\in \cP$ then, by 
Remark~\ref{ref-def-adm}\eqref{ref-def-adm2}, $v_i\Vert v_j$. 

\item $\ell_i,\ell_j\in \cL_1$; considering 
any $\ell_k\in \cL_2$ and using the previous
case, one has $v_i\Vert v_k\Vert v_j$. The same argument
works for $\ell_i,\ell_j\in \cL_2$. In
particular, $v_i\Vert v_j$ if $\ell_i,\ell_j\in \cL_1\cup\cL_2$.

\item $\ell_i\in\cL_0$ and $P\in\cP$ such that $P\cap\cL_0=\{\ell_i\}$. 
Since all the vectors associated with $P$ but one are proportional, then
this must also be the case for $v_i$.
\end{itemize}
All the vectors (but at most one) are proportional. To conclude
we apply condition~\eqref{def-adm-3} of Definition~\ref{def-adm}.

\item\label{exam-adm-ceva}
Ceva's line combinatorics is 3-admissible.
\begin{proof}
Ceva's line combinatorics is given by the following realization:
\begin{figure}[ht]
\includegraphics[scale=.25]{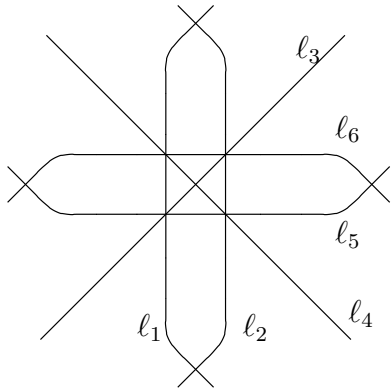}
\begin{picture}(0,0)
\put(-60,20){$\ell_2$}
\put(-20,25){$\ell_4$}
\put(-25,55){$\ell_5$}
\put(-25,95){$\ell_6$}
\put(-40,125){$\ell_3$}
\put(-100,20){$\ell_1$}
\end{picture}
\caption{Ceva's line combinatorics}
\end{figure}
$$\cL_{\text{\sc ceva}}:=\{1,2,3,4,5,6\}$$
$$\cP_{\text{\sc ceva}}:=\{\{1,2\},\{3,4\},\{5,6\},\{1,3,5\},\{1,4,6\},\{2,3,6\},\{2,4,5\}\}.$$
For example, the following is a 3-admissible set of vectors for Ceva:
$$\left\{ v_1=v_2=(1,0), v_3=v_4=(0,1), v_5=v_6=(-1,-1) \right\}$$
\end{proof}
\item\label{exam-adm-mcl}
MacLane's line combinatorics $\mlc$ is not 3-admissible.
\begin{proof}
We will use the combinatorics given in Figure~\ref{mclane}. Let us assume that MacLane 
is 3-admissible, then one has a list of non-zero vectors $v_0,v_1,...,v_7$ associated 
with each line. We will first see that $v_0$ and $v_1$ cannot be proportional. If they 
were ($v_0\Vert v_1$), using $\{0,1,2\}$, $\{1,6\}$ and $\{2,3\}$ one would have that 
$v_0\Vert v_2\Vert v_6\Vert v_3$ and finally, using $\{3,6,7\}$, one obtains 
$v_0\Vert v_7$, and hence, by 
$\{1,4,7\}$ and $\{4,5\}$, all vectors are proportional, which contradicts 
condition~\eqref{def-adm-1} of Definition~\ref{def-adm}.

Therefore, $v_0$ and $v_1$ are linearly independent. After a change of basis,
one can assume that $v_0 = (\alpha,0)$, $v_1=(0,\beta)$, $\alpha,\beta\in\ZZ\setminus\{0\}$,
and therefore $v_2=(-\alpha,-\beta)$. We will 
briefly describe the conditions and how they affect the vectors until a contradiction 
is reached:
$$
\rightmap{{\tiny{\array{c}\{0,7\}\\ \{1,6\}\endarray}}} 
\begin{cases}
v_7=(\gamma ,0) \\
v_6=(0,\delta)
\end{cases}
\rightmap{{\tiny{\array{c}\{3,6,7\}\\ \{0,5,6\}\endarray}}} 
\begin{cases}
v_3 = (-\gamma ,-\delta)\\
v_5 = (-\alpha,-\delta)\\
\end{cases}
\rightmap{\{1,3,5\}}
\begin{cases}
\alpha =-\gamma\\
\beta = 2\delta.
\end{cases}
$$

Hence $v_3$  does not satisfy
the condition~\eqref{def-adm-1b} of Definition~\ref{def-adm} with $v_2=(\gamma,-2\delta)$ 
on the double point $\{2,3\}$.
\end{proof}
\end{enumerate}
\end{examples}

\begin{definition}
\label{def-pt-adm}
A combinatorics $\scc$ with only double and triple points 
is \emph{pointwise 3-admissible} if
the only 3-admissible subcombinatorics of $\scc$ are isomorphic to $\cM_3$.
\end{definition}

\begin{proposition}
\label{propos-permut}
If $\scc$ is a pointwise 3-admissible combinatorics 
then any $\psi\in \Aut^1(H_\scc)$ induces a permutation 
$\psi_3$ of $\cP_3$.
\end{proposition}

\begin{proof}
Let $\cL:=\{\ell_0,...,\ell_r\}$ denote the set of lines of $\scc$ and let
$\cP_3\subset \cP$ denote the set of triple points of $\scc$. We will first prove that 
any isomorphism $\psi\in \Aut^1(H_\scc)$ induces a map $\psi_3:\cP_3\to \cP_3$. 
Consider a triple point $P:=\{\ell_i,\ell_j,\ell_k\}\in \cP_3$;
then $\Adm_\psi(P)$ is an admissible subcombinatorics
of $\scc$ (Example~\ref{exam-adm}\eqref{exam-adm-adm})
and defines a triple
point. The map $\psi_3:\cP_3\to \cP_3$ given by $\psi_3(P):=\Adm_\psi(P)$
is defined.

We will next prove that such a map is indeed injective, and hence a permutation (in order
for this to make sense, one has to assume that $\# \cP_3>1$ and hence $r\geq 4$). 
Assume $\psi_3$ is not injective, and let 
$\psi_3(P_1)=\psi_3(P_2)=Q=\{\ell_u,\ell_v,\ell_w\}$. One has to consider
two different cases depending on whether or not $P_1$ and $P_2$ 
share a line:
\begin{enumerate}
\item
If $P_1,P_2$ do not share a line, i.e, $P_1:=\{\ell_{i_1},\ell_{j_1},\ell_{k_1}\}$ and 
$P_2:=\{\ell_{i_2},\ell_{j_2},\ell_{k_2}\}$, where all the subscripts are pairwise different. 
By reordering the columns, let us write $Q=\{\ell_0,\ell_1,\ell_2\}$. Let
$A^\psi_{P_1,P_2}$ be the submatrix of $A^\psi$ 
corresponding to the rows $\{i_1,j_1,k_1,i_2,j_2,k_2\}$.
Using \eqref{eq-unos}:
\begin{equation}
\label{eq-matrix1}
A^\psi_{P_1,P_2}:=
\begin{pmatrix}
a^{i_1}_0 & a^{i_1}_1 & a^{i_1}_2 & a^{i_1}_3 & \dots & a^{i_1}_{r} \\
a^{j_1}_0 & a^{j_1}_1 & a^{j_1}_2 & a^{i_1}_3 & \dots & a^{i_1}_{r} \\
a^{k_1}_0 & a^{k_1}_1 & a^{k_1}_2 & a^{i_1}_3 & \dots & a^{i_1}_{r} \\
a^{i_2}_0 & a^{i_2}_1 & a^{i_2}_2 & a^{i_2}_3 & \dots & a^{i_2}_{r} \\
a^{j_2}_0 & a^{j_2}_1 & a^{j_2}_2 & a^{i_2}_3 & \dots & a^{i_2}_{r} \\
a^{k_2}_0 & a^{k_2}_1 & a^{k_2}_2 & a^{i_2}_3 & \dots & a^{i_2}_{r} \\
\end{pmatrix},
\end{equation}
where 
$a^{i_\bullet}_0+a^{i_\bullet}_1+a^{i_\bullet}_2
=a^{j_\bullet}_0+a^{j_\bullet}_1+a^{j_\bullet}_2
=a^{k_\bullet}_0+a^{k_\bullet}_1+a^{k_\bullet}_2$, $\bullet=1,2$. 
The sublattice $K$ of $\Sigma_5$ generated by
the columns $\!\!\pmod{\uno_6}$ should have maximal rank equal
to $5$. Note that $\rank(K)$ equals the rank of the matrix
$\overline A^\psi_{P_1,P_2}$ obtained by subtracting the last row
from the first ones, forgetting the
last row and replacing the first column by
the sum of the first three:
\begin{equation}
\label{eq-matrix1b}
\overline A^\psi_{P_1,P_2}=
\begin{pmatrix}
b^{i_1}_0 & b^{i_1}_1 & b_2 & b_3 & \dots & b_{r} \\
b^{j_1}_0 & b^{j_1}_1 & b_2 & b_3 & \dots & b_{r} \\
b^{k_1}_0 & b^{k_1}_1 & b_2 & b_3 & \dots & b_{r} \\
b^{i_2}_0 & b^{i_2}_1 & 0   & 0   & \dots & 0 \\
b^{j_2}_0 & b^{j_2}_1 & 0   & 0   & \dots & 0 \\
\end{pmatrix},
\end{equation}
which does not have rank $5$.
\item
If $P_1,P_2$ share a line, say $P_1:=\{i,j_1,k_1\}$ and $P_1:=\{i,j_2,k_2\}$.
Then, analogously to the previous case, one obtains a similar matrix to
(\ref{eq-matrix1}) but where the rows $i_1$ and $i_2$ are identified,
and we proceed in a similar way.
\end{enumerate}
\end{proof}

\begin{definition}
Three triple points $P,Q,R\in \cP$ of a line combinatorics $(\cL,\cP)$ are
said to be \emph{in a triangle} if $P\cap Q=\{\ell_{1}\}$, 
$P\cap R=\{\ell_{2}\}$ and $Q\cap R=\{\ell_{3}\}$ are pairwise different.
\end{definition}

\begin{proposition}
\label{propos-c3}
For any $\psi\in \Aut^1(H_\scc)$, $\scc$ pointwise $3$-admissible,
$\psi_3$  satisfies the following Triangle Property:
$\psi_3:\cP_3\to \cP_3$ preserves triangles, that is,
if $P_1,P_2,P_3\in \cP_3$ are in a triangle, then
$\psi_3(P_1),\psi_3(P_2),\psi_3(P_3)$ are also in a triangle.

\end{proposition}

\begin{proof}
Let $P_1,P_2,P_3\in \cP_3$ be three triple points in a triangle, say 
$P_1:=\{\ell_i,\ell_j,\ell_k\}$,
$P_2:=\{\ell_k,\ell_l,\ell_m\}$, $P_3:=\{\ell_m,\ell_n,\ell_i\}$.
Let us assume that $\psi_3(P_1),\psi_3(P_2),\psi_3(P_3)$ are not in a triangle.
One has two possibilities, either two of them do not share a line or three of them share 
a line. 
\begin{enumerate}
\item Two of them, say $\psi_3(P_1),\psi_3(P_2)$ do not share a line.
After reordering, we can suppose that $\psi_3(P_1)=\{\ell_0,\ell_1,\ell_2\}$ and 
$\psi_3(P_2)=\{\ell_3,\ell_4,\ell_5\}$. For $\psi_3(P_3)$
there are several possibilities but we may assume
$\psi_3(P_3)\subset\{\ell_0,\ell_3,\ell_6,\ell_7,\ell_8\}$.
Consider 
the submatrix $A^\psi_{P_1,P_2,P_3}$ of $A^\psi$ given by the rows of $P_1,P_2,P_3$. 
Applying~\eqref{eq-unos} successively to the rows
defined by $P_1,P_2,P_3$ one has:
\begin{equation}
\label{eq-matrix2b}
\setcounter{MaxMatrixCols}{12}
A^\psi_{P_1,P_2,P_3}:=
\begin{pmatrix}
a^i_0& a^i_1& a^i_2& a_3& a_4& a_5& a_6& a_7& a_8& a_9& \dots & a_{r}\\
a^j_0& a^j_1& a^j_2& a_3& a_4& a_5& a_6& a_7& a_8& a_9& \dots & a_{r}\\
a_0& a_1& a_2& a_3& a_4& a_5& a_6& a_7& a_8& a_9& \dots & a_{r}\\
a_0& a_1& a_2& a^l_3& a^l_4& a^l_5& a_6& a_7& a_8& a_9& \dots & a_{r}\\
a_0& a_1& a_2& a^m_3& a_4& a_5& a_6& a_7& a_8& a_9& \dots & a_{r}\\
a^n_0& a_1& a_2& a^n_3& a_4& a_5& a^n_6& a^n_7& a^n_8& a_9& \dots & a_{r}\\
\end{pmatrix};
\end{equation}
moreover, $a^i_1=a_1$, $a^i_2=a_2$,
\begin{alignat*}{3}
a^i_0+a_1+a_2&=a^j_0+a^j_1+a^j_2&=a_0+a_1+a_2& (\Rightarrow a^i_0=a_0),\\
a_3+a_4+a_5&=a^l_3+a^l_4+a^l_5&=a^m_3+a_4+a_5& (\Rightarrow a^m_3=a_3)
\end{alignat*}
and
$$a_0+a_3+a_6+a_7+a_8=a^n_0+a_3^n+a^n_6+a^n_7+a^n_8.$$
As in the proof of Proposition~\ref{propos-permut}
we need $\rank(\overline A^\psi_{P_1,P_2,P_3})=5$,
where $\overline A^\psi_{P_1,P_2,P_3}$ is
the matrix obtained from $A^\psi_{P_1,P_2,P_3}$
by subtracting the last row
from the first ones and forgetting the last row.
We obtain:
\begin{equation}
\label{eq-matrix2c}
\setcounter{MaxMatrixCols}{12}
\bar A^\psi_{P_1,P_2,P_3}:=
\begin{pmatrix}
b_0  & 0    & 0    & b_3  & 0   &0    &b_6&b_7&b_8& 0& \dots & 0\\
b^j_0& b^j_1& b^j_2& b_3  & 0   & 0   &b_6&b_7&b_8& 0& \dots & 0\\
b_0  & 0    & 0    & b_3  & 0   &0    &b_6&b_7&b_8& 0& \dots & 0\\
b_0  & 0    & 0    & b^l_3&b^l_4&b^l_5&b_6&b_7&b_8& 0& \dots & 0\\
b_0  & 0    & 0    & b_3  &0    &0    &b_6&b_7&b_8& 0& \dots & 0\\
\end{pmatrix}
\end{equation}
which cannot have rank $5$.
\item If $\psi_3(P_1),\psi_3(P_2)$ and $\psi_3(P_3)$ have a common line, say
$\ell_0$, we follow the same strategy and obtain the desired result.
\end{enumerate}
\end{proof}

Our main goal is to check if $\psi_3$ is induced by an element
of $\Aut(\scc)$. The next example shows that we need \emph{enough}
triangles.

\begin{example}
Note that Proposition~\ref{propos-permut} does not automatically ensure that
in general an automorphism of the combinatorics is produced. 
For instance, consider the combinatorics 
$\scc$ given by the lines $\{0,1,\ldots,6\}$, and the following triple points 
$\{0,1,2\}$, $\{2,3,4\}$, and $\{4,5,6\}$ (the remaining intersections are double points). 
It is easy to see that such a combinatorics is pointwise 3-admissible.
Let $\psi:H_\scc\to H_\scc$ be given by the following matrix:
$$
A^\psi:=\left(
\begin{array}{ccccccc}
1&0 & 0 & 0 & 0 & 0 & 0\\
0&1 & 0 & 0 & 0 & 0 & 0 \\ 
0&0 & 1 & 0 & 0 & 0 & 0 \\ 
0&0 & 1 & 0 & 1 & 0 & -1 \\ 
0&0 & 1 & 0 & 0 & 1 & -1 \\ 
0&0 & 0 & 1 & 0 & 1 & -1 \\ 
0&0 & 0 & 0 & 1 & 1 & -1
\end{array} \right),
$$
satisfying the admissibility relations. This induces the following maps:
$$\array{cccc}
\array{ccc}
\cP_3 & \rightmap{\psi_3} & \cP_3\\
\{0,1,2\}&\mapsto &\{0,1,2\}\\
\{2,3,4\}&\mapsto &\{4,5,6\}\\
\{4,5,6\}&\mapsto &\{2,3,4\}
\endarray
&\quad & \quad &
\array{ccc}
\gr^1M_\scc & \rightmap{\wedge^2 \psi} & \gr^1M_\scc\\
x_{0,1}&\mapsto &x_{0,1}\\
x_{2,3}&\mapsto &x_{4,5}\\
x_{4,5}&\mapsto &x_{2,3},
\endarray
\endarray$$
where 
$\gr^1 M_\scc\cong\langle x_{0,1} \rangle \ZZ \oplus \langle x_{2,3} \rangle \ZZ \oplus
\langle x_{4,5} \rangle \ZZ$.

However, the given permutation is not induced by an 
automorphism of the combinatorics, because the point $\{2,3,4\}$ (which is 
the only one that shares a line with the other two) is not fixed.
\end{example}

We want to apply the previous results to $\rybc$. First we
will check that $\rybc$ is pointwise 3-admissible.

\begin{lemma}
\label{lem-ryb-mcl}
An admissible subcombinatorics of $\rybc$ cannot have lines in both 
$\cR_1:=\{3,4,5,6,7\}$ and $\cR_2:=\{8,9,10,11,12\}$.
\end{lemma}
\begin{proof}
Any subcombinatorics of $\rybc$ having lines in both $\cR_1$ and $\cR_2$
verifies the conditions of Example~\ref{exam-adm}(\ref{exam-adm-general}).
\end{proof}

\begin{lemma}
\label{lem-mac-no-adm}
$\mlc$ is pointwise 3-admissible.
\end{lemma}
\begin{proof}
It is not difficult to prove that any combinatorics (with only double and triple
points) of less than 6 lines (other than $\cM_3$) is not $3$-admissible. In
Example~\ref{exam-adm}\eqref{exam-adm-mcl} it is shown for the whole combinatorics of 
eight lines. Let us check the remaining cases, that is, six and seven lines. 
Note that, up to a 
combinatorics automorphism, there is a unique way to remove one line. There are,
however, two possible ways to remove two lines, depending on whether they intersect
or not at a triple point. Therefore we only have to check the following cases:
\begin{enumerate}
\item For 7 lines \\
$\{0,1,2,3,4,5,6\}$: 
$\rightmap{\{3,6\}} v_3\Vert v_6 \rightmap{\{1,6\},\{2,3\}} 
v_3\Vert v_1\Vert v_2 \rightmap{\{1,4\}} v_3\Vert v_4 \rightmap{\{4,5\}} v_3\Vert v_5$.
\item For 6 lines, removing two lines intersecting at a triple point \\
$\{0,2,3,4,5,6\}$:
$\rightmap{\{3,5\}} v_5\Vert v_3 \rightmap{\{2,3\},\{3,6\},\{4,5\}} 
v_5\Vert v_2\Vert v_6\Vert v_4$.
\item For 6 lines, removing two lines intersecting at a double point \\
$\{0,2,3,4,5,7\}$:
$\rightmap{\{2,3\},\{2,4\}} v_2\Vert v_3\Vert v_4 \rightmap{\{3,7\},\{4,5\}} 
v_2\Vert v_5\Vert v_7$.
\end{enumerate}
\end{proof}

\begin{proposition}
\label{props-ryb-no-adm}
$\rybc$ is pointwise 3-admissible.
\end{proposition}
\begin{proof}
An immediate consequence of Lemmas~\ref{lem-ryb-mcl} and \ref{lem-mac-no-adm}.
\end{proof}

\begin{remarks}
\label{rem-ryb-triang}
In order to prove that for any $\psi\in\Aut^1(H_\rybc)$, $\psi_3$
comes from an element of $\Aut(\rybc)$ we need to know more 
combinatorial properties
of $\rybc$ and $\mlc$.
\begin{enumerate}
\item\label{rem-ryb-triang36}
The triple point $\{0,1,2\}\in \cP_3$ in $\rybc$ is the only one that belongs to 36 triangles.
\item\label{rem-ryb-triang-triple}
Any triple point $P\in \cP_3$ in $\rybc$ except for $\{3,6,7\}$ and $\{8,11,12\}$ 
satisfies that $\{0,1,2\}$ and $P$ are in a triangle.
\item\label{rem-ryb-triang-2triple}
Any two triple points in $\mlc$ sharing a line belong to a triangle.
\item\label{rem-ryb-triang-3triple}
For any three triple points $P_1,P_2,P_3$ in $\mlc$, there exists another
triple point $Q$ such that $Q,P_i,P_j$ belong to a triangle ($i,j\in \{1,2,3\}$).
\end{enumerate}
\end{remarks}

\begin{proposition}
\label{cor-mclane} Let $\psi\in\Aut^1(\rybc)$.
\begin{enumerate}[\rm(1)]

\item\label{cor-mclane1} $\psi_3(\{0,1,2\})=\{0,1,2\}$ 
\item\label{cor-mclane2} $\psi_3$ either preserves (resp. exchanges) the triple points
of $\cR_1$ and $\cR_2$ in $\rybc$ inducing an automorphism
\item\label{cor-mclane3} The action of $\psi_3$ on $\cR_1$ and $\cR_2$
comes from an automorphism (resp. isomorphism) of their combinatorics.
\end{enumerate}
\end{proposition}

\begin{proof}
Part \eqref{cor-mclane1} is true by Propositions~\ref{propos-permut}-\ref{propos-c3} and 
Remark~\ref{rem-ryb-triang}(\ref{rem-ryb-triang36}).
By Remark~\ref{rem-ryb-triang}(\ref{rem-ryb-triang-triple}), the points $\{3,6,7\}$ 
and $\{8,11,12\}$ are either preserved or exchanged. In order to prove~\eqref{cor-mclane2} 
and~\eqref{cor-mclane3}, we may suppose that $\psi_3(\{3,6,7\})=\{3,6,7\}$.
Recall that the subcombinatorics defined by $\cR_0\cup\cR_1$
is isomorphic to $\mlc$.
Since triangles are preserved by $\psi_3$ (Proposition~\ref{propos-c3}), 
according to Remark~\ref{rem-ryb-triang}(\ref{rem-ryb-triang-2triple}) the images of 
any two triple points in $\cR_0\cup\cR_1$ sharing a line also share a line.
This implies, using Remark~\ref{rem-ryb-triang}(\ref{rem-ryb-triang-3triple}) again, that
the image of any three triple points on $\cR_0\cup\cR_1$ sharing a line are also three points 
sharing a line. Since any line in $\cR_1$ or $\cR_2$ has at least three triple points, we 
conclude~\eqref{cor-mclane2} and~\eqref{cor-mclane3}.
\end{proof}

\begin{proposition}
\label{propos-combinat}
Let $\psi\in \Aut^1(H_\rybc)$; then
$\psi_3:\cP_3\to \cP_3$ is induced by an automorphism of $\rybc$.
\end{proposition}

\begin{proof}
We will use Proposition~\ref{cor-mclane} repeatedly.
We can compose $\psi$ with an element of $\Aut(\rybc)$
in order to have $\psi_3$ preserve the triple
points in $\cR_i$. Recall that $\psi_3|_{\cR_0\cup\cR_i}$
comes from an automorphism $\varphi_i$ of $\mlc$ which
respects $\{0,1,2\}$. Composing again with an
element of $\Aut(\rybc)$ we may suppose that
$\varphi_1$ is the identity on $\{0,1,2\}$.
It is enough to prove that it is also the case for $\varphi_2$.
If it is not the case, we may assume (by conjugation
with an element of $\Aut(\rybc)$) that
$\varphi_2(0)=1$.
There are two possibilities to be checked, depending on whether
$9$ and $10$ are fixed or permuted. In both cases, the triple points
$\{0,1,2\}$, $\{3,6,7\}$ and $\{8,11,12\}$ are fixed. Using
the arguments in the proofs of Propositions~\ref{propos-permut} and~\ref{propos-c3},
one can obtain the induced matrices 
$A^\psi$ with all the admissibility relations, which, modulo $\uno_{13}$ do not have a
maximal rank in either case.
\end{proof}

\begin{proposition}
\label{thm-main-comb}
$\rybc$ is homologically rigid.
\end{proposition}

\begin{proof}
Let $\psi\in \Aut^1(H_\rybc)$; by Propositions~\ref{propos-combinat}
and~\ref{propos-c3}, 
$\psi_3$ comes from an automorphism of the combinatorics. 
Composing with the inverse of such an automorphism, we may suppose
that $\psi_3=1_{\cP_3}$. It is hence enough 
to prove that any isomorphism $\psi\in \Aut^1(H_\rybc)$ that induces the identity on 
$\psi_3$ is just $\pm 1_{H_\rybc}$. 
From the definition of $\Adm_\psi(P)$, $P\in\cP_3$, we deduce
the following. Let us fix the $j^{th}$-column; all the entries
in this column corresponding to $P\in\cP_3$ such that $j\notin P$
are equal. We deduce from this that we can choose $A^\psi$
such that all the elements outside the diagonal are constant
in their column. Adding multiples of $\uno_{13}$, we obtain
that  $A^\psi$ can be chosen to be diagonal. We also know
that for each $P\in\cP_3$, the diagonal terms corresponding
to $\cP$ are equal and since any two elements can be joined
by a \emph{chain} of triple points, we deduce that all the
diagonal terms are equal. Since $\psi$ is an automorphism,
they are equal to $\pm1$ and then  $\psi=\pm 1_{H_\rybc}$.
\end{proof}

Therefore we can prove the main result.

\begin{theorem}
The fundamental groups of the two complex realizations of Rybnikov's combinatorics 
are not isomorphic.
\end{theorem}

\begin{proof}
Let $G_{+}$ and $G_{-}$ be the fundamental groups of $R_{\gw,\gw}$ and 
$R_{\bar \gw,\gw}$ respectively. Any isomorphism $\tilde\psi:G_{+}\to G_{-}$ will produce 
an automorphism $\psi:H_\rybc\equiv G_{+}/G'_{+}\to G_{-}/G'_{-}\equiv H_\rybc$, that is,
we can consider $\psi\in \Aut^1(\rybc)$. By Theorem~\ref{thm-main-comb}, $\psi$ induces an
automorphism of $\rybc$. Since the identifications $H_\rybc\equiv G_{\pm}/G'_{\pm}$ 
are made up to the action of $\Aut(\rybc)$, we may assume that $\psi$ induces 
$\pm 1_{H_\rybc}$.
Moreover, eventually exchanging $R_{\gw,\gw}$ (resp. $R_{\bar \gw,\gw}$) for
$R_{\bar\gw,\bar\gw}$ (resp. $R_{\gw,\bar\gw}$), see Example~\ref{exam-ryb-comb} by
means of the automorphism given by complex conjugation, we may assume that $\psi=1_{H_\rybc}$
(Proposition~\ref{thm-main-comb}). Therefore $\tilde\psi$ is a homologically trivial 
isomorphism between $G_{+}$ and $G_{-}$, something which is ruled out by
Theorem~\ref{thm-main-ryb}.
\end{proof}

\def\cprime{$'$}
\providecommand{\bysame}{\leavevmode\hbox to3em{\hrulefill}\thinspace}
\providecommand{\MR}{\relax\ifhmode\unskip\space\fi MR }
\providecommand{\MRhref}[2]{%
  \href{http://www.ams.org/mathscinet-getitem?mr=#1}{#2}
}
\providecommand{\href}[2]{#2}

\end{document}